# ON THE PAPER "WEAK CONVERGENCE OF SOME CLASSES OF MARTINGALES WITH JUMPS"

By Yoichi Nishiyama

*Institute of Statistical Mathematics*

This note extends some results of Nishiyama [*Ann. Probab.* **28** (2000) 685–712]. A maximal inequality for stochastic integrals with respect to integer-valued random measures which may have infinitely many jumps on compact time intervals is given. By using it, a tightness criterion is obtained; if the so-called *quadratic modulus* is bounded in probability and if a certain entropy condition on the parameter space is satisfied, then the tightness follows. Our approach is based on the entropy techniques developed in the modern theory of empirical processes.

**1. Introduction.** This note extends a maximal inequality of [3], and develops the entropy methods for martingales which have been studied systematically in [4]. Let $(E, \mathcal{E})$ be a Blackwell space. For every $n \in \mathbb{N}$, let $\mu^n$ be an integer-valued random measure on $\mathbb{R}_+ \times E$ defined on a stochastic basis $\mathbf{B}^n = (\Omega^n, \mathcal{F}^n, \mathbf{F}^n = (\mathcal{F}^n_t)_{t \in \mathbb{R}_+}, P^n)$, and let $\nu^n$ be the predictable compensator of $\mu^n$. Let $\mathcal{W}^n = \{W^{n,\psi} : \psi \in \Psi\}$ be a class of predictable functions on $\Omega^n \times \mathbb{R}_+ \times E$, indexed by an arbitrary set $\Psi$. Let $\tau^n$ be a finite stopping time. We treat the sequence of processes $(t, \psi) \rightsquigarrow X^{n,\psi}_t$ given by

$$X^{n,\psi}_t = W^{n,\psi} * (\mu^n - \nu^n)_t \qquad \forall t \in \mathbb{R}_+, \ \forall \psi \in \Psi.$$

[2] and [3] studied the weak convergence of the sequence as $n \to \infty$, and the latter showed that, if the so-called *quadratic modulus* is bounded in probability and if a certain integrability condition for "partitioning entropy" is satisfied, then the tightness in $\ell^\infty([0, t_0] \times \Psi)$, where $t_0 = \tau^n$ is a constant, of the sequence is implied. However, [3] assumed the following:

*Case* $\mathrm{A}^n$. The process $t \rightsquigarrow \overline{W}^n * \nu^n_t$ is locally integrable and $\nu^n([0, \tau^n] \times E) < \infty$ almost surely.

Here, "$\overline{W}^n = \sup_{\psi \in \Psi} |W^{n,\psi}|$" [the meaning of the quotation mark is that, precisely speaking, we have to take a "predictable envelope"]. The assumption that $\nu^n([0, \tau^n] \times E) < \infty$ a.s. implies that only finitely many jumps of the process may occur. Although such a case already serves a lot of applications, the above situation is far from the general theory of local martingales. The main contribution of this note is to replace the above assumption by the following two conditions:

*Case* B$^n$. $t \rightsquigarrow (|\overline{W}^n|^2 \wedge \overline{W}^n) * \nu_t^n$ is locally integrable; $\Psi$ is countable and is "asymptotically separated" by a series of finite partitions.

The definition of the notion "asymptotically separate" will be given in Section 2.

Our result is of interest by itself. We refer to [3] for the further discussions on the background of our results and the history of related works. In addition to the theoretical interest, the processes with infinitely many jumps have recently been important in applications, for example, in the context of mathematical finance. Although this short note is not a place to present technical examples, our result would hopefully yield some new applications, especially semi- and nonparametric statistical inferences. Actually, an application to Lévy processes has been established in [5].

In Section 2 we give some additions to [3] for which we follow all definitions and notation. A maximal inequality, which has the same form as [3], is given in the case where infinitely many jumps may occur. The change of the proof is just one point, so we try to reach there in as few pages as possible, and to explain the difference clearly. By using the inequality, we give a sufficient condition for the processes to take values in $\ell^\infty$-spaces. We state a weak convergence theorem which is an immediate consequence of those results. Proofs are given in the Appendix.

## 2. Results.   Let us begin with preparing three definitions.

DEFINITION 2.1.   Let $(\mathcal{X}, \mathcal{A}, \lambda)$ be a $\sigma$-finite measure space. For a given mapping $Z : \mathcal{X} \to \mathbb{R} \cup \{\infty\}$, we denote by $[Z]_{\mathcal{A}, \lambda}$ any $\mathcal{A}$-measurable function $U : \mathcal{X} \to \mathbb{R} \cup \{\infty\}$ such that: (i) $U \geq Z$ holds identically; (ii) $\tilde{U} \geq U$ holds $\lambda$-almost everywhere, for every $\mathcal{A}$-measurable function $\tilde{U}$ such that $\tilde{U} \geq Z$ holds $\lambda$-almost everywhere.

The existence of such a random variable $[Z]_{\mathcal{A}, \lambda}$ and its uniqueness up to a $\lambda$-negligible set follow from Lemma 1.2.1 of [6].

DEFINITION 2.2.   Let $\Psi$ be an arbitrary set. $\Pi = \{\Pi(\varepsilon)\}_{\varepsilon \in (0, \Delta_\Pi]}$, where $\Delta_\Pi \in (0, \infty) \cap \mathbb{Q}$, is called a decreasing series of finite partitions (abb. DFP) [resp., nested series of finite partitions (abb. NFP)] of $\Psi$ if it satisfies the following (i), (ii) and (iii) [resp., (i), (ii) and (iii$'$)]: (i) each $\Pi(\varepsilon) = \{\Psi(\varepsilon; k) : 1 \leq$ $k \leq N_\Pi(\varepsilon)\}$ is a finite partition of $\Psi$, that is, $\Psi = \bigcup_{k=1}^{N_\Pi(\varepsilon)} \Psi(\varepsilon; k)$; (ii) $N_\Pi(\Delta_\Pi) = 1$ and $\lim_{\varepsilon \downarrow 0} N_\Pi(\varepsilon) = \infty$; (iii) $N_\Pi(\varepsilon) \geq N_\Pi(\varepsilon')$ whenever $\varepsilon \leq \varepsilon'$; (iii') $\Pi(\varepsilon) \supset \Pi(\varepsilon')$ whenever $\varepsilon \leq \varepsilon'$.

Notice that any NFP is a DFP.

DEFINITION 2.3.   We say a DFP $\Pi = \{\Pi(\varepsilon)\}_{\varepsilon \in (0, \Delta_\Pi]}$ of $\Psi$ asymptotically separates $\Psi$ if for any finite subset $F \subset \Psi$ there exists $\varepsilon_F$ such that for every $\varepsilon \in (0, \varepsilon_F]$ each partitioning set $\Psi(\varepsilon; k)$ of the partition $\Pi(\varepsilon)$ contains at most one point of $F$.

This is not a strong requirement. In fact, consider the case where $\Psi$ is totally bounded with respect to a metric $\rho$. When each $\Pi(\varepsilon)$ is a partition generated by $\varepsilon$-balls which cover $\Psi$, then $\Pi = \{\Pi(\varepsilon)\}_{\varepsilon \in (0, \Delta_\Pi]}$ asymptotically separates $\Psi$.

Let us now turn to the context of integer-valued random measure. Let $(E, \mathcal{E})$ be a Blackwell space. Let $\mu$ be an integer-valued random measure on $\mathbb{R}_+ \times E$ defined on a stochastic basis $\mathbf{B} = (\Omega, \mathcal{F}, \mathbf{F} = (\mathcal{F}_t)_{t \in \mathbb{R}_+}, P)$, and $\nu$ a "good" version of the predictable compensator of $\mu$. Let $\tau$ be a finite stopping time. We put $\tilde{\Omega} = \Omega \times \mathbb{R}_+ \times E$ and $\tilde{\mathcal{P}} = \mathcal{P} \otimes \mathcal{E}$, where $\mathcal{P}$ is the predictable $\sigma$-field. Let $\mathcal{W} = \{W^\psi : \psi \in \Psi\}$ be a family of predictable functions on $\tilde{\Omega}$ indexed by $\Psi$. We introduce the Doléans measure $M_\nu^P$ on $(\tilde{\Omega}, \tilde{\mathcal{P}})$, which is $\tilde{\mathcal{P}}$-$\sigma$-finite, given by

$$M_\nu^P(d\omega, dt, dx) = P(d\omega)\nu(\omega; dt, dx).$$

(See Section II.1 of [1] for the theory of random measures.)

Let us recall the definitions of the *predictable envelope* $\overline{W}$ and the *quadratic* $\Pi$-*modulus* $\|\mathcal{W}\|_\Pi$ given by [3].

DEFINITION 2.4.   The predictable envelope $\overline{W}$ of $\mathcal{W} = \{W^\psi : \psi \in \Psi\}$ is defined by

$$\overline{W} = \left[\sup_{\psi \in \Psi} |W^\psi|\right]_{\tilde{\mathcal{P}}, M_\nu^P}.$$

For a given DFP $\Pi$ of $\Psi$, the quadratic $\Pi$-modulus $\|\mathcal{W}\|_\Pi$ of $\mathcal{W} = \{W^\psi : \psi \in \Psi\}$ is defined as the $\mathbb{R}_+ \cup \{\infty\}$-valued predictable process $t \rightsquigarrow \|\mathcal{W}\|_{\Pi, t}$ given by

$$\|\mathcal{W}\|_{\Pi, t} = \sup_{\varepsilon \in (0, \Delta_\Pi] \cap \mathbb{Q}} \max_{1 \leq k \leq N_\Pi(\varepsilon)} \frac{\sqrt{|\Delta_W(\Psi(\varepsilon; k))|^2 * \nu_t}}{\varepsilon} \qquad \forall t \in \mathbb{R}_+,$$

where

$$\Delta_W(\Psi') = \left[\sup_{\psi, \phi \in \Psi'} |W^\psi - W^\phi|\right]_{\tilde{\mathcal{P}}, M_\nu^P} \qquad \forall \Psi' \subset \Psi.$$

We will consider the two cases:

*Case* A. The process $t \rightsquigarrow \overline{W} * \nu_t$ is locally integrable and $\nu([0, \tau] \times E) < \infty$ almost surely.

*Case* B. The process $t \rightsquigarrow (\overline{W}^2 \wedge \overline{W}) * \nu_t$ is locally integrable, $\Psi$ is countable, and $\Pi$ is a DFP which asymptotically separates $\Psi$.

In any case, we define the following:

$$(1) \qquad X_t^\psi = W^\psi * (\mu - \nu)_t \qquad \qquad \forall \psi \in \Psi;$$

$$(2) \qquad X_t^{a,\psi} = W^\psi 1_{\{\overline{W} \leq a\}} * (\mu - \nu)_t \qquad \forall \psi \in \Psi \ \forall a > 0;$$

$$(3) \qquad \check{X}_t^{a,\psi} = W^\psi 1_{\{\overline{W} > a\}} * (\mu - \nu)_t \qquad \forall \psi \in \Psi \ \forall a > 0.$$

Our main interest is the process $t \rightsquigarrow X_t^\psi$. In both cases, the process $t \rightsquigarrow X_t^\psi$ is a local martingale and the process $t \rightsquigarrow X_t^{a,\psi}$ is a locally square-integrable martingale. In Case A they have finite variation, while in Case B they may not. In both cases the process $t \rightsquigarrow \check{X}_t^{a,\psi}$ is a local martingale which has finite variation. (See Proposition II.1.28 and Theorem II.1.33 of [1].) Notice that, in Case B, the equality like $W * (\mu - \nu)_t = W * \mu_t - W * \nu_t$ may *not* hold, and that the processes $X^\psi$ and $X^{a,\psi}$ are defined for all $\psi \in \Psi$, only *almost surely*; see [3].

The following theorem gives some maximal inequalities for these processes in terms of $\|\mathcal{W}\|_\Pi$. Here, the notation "$\lesssim$" means that the left-hand side is not bigger than the right-hand side up to a multiplicative universal constant.

**THEOREM 2.5.** *The following* (i) *and* (ii) *hold not only in Case* A *but also in Case* B.

(i) *For given NFP* $\Pi$ *of* $\Psi$ *and any constants* $\delta \in (0, \Delta_\Pi)$ *and* $K > 0$,

$$E^* \sup_{t \in [0,\tau]} \max_{1 \leq k \leq N_\Pi(\delta)} \sup_{\psi, \phi \in \Psi(\delta;k)} |X_t^{a,\psi} - X_t^{a,\phi}| 1_{\{\|\mathcal{W}\|_{\Pi,\tau} \leq K\}}$$

$$\lesssim K \int_0^\delta \sqrt{\log(1 + N_\Pi(\varepsilon))} \, d\varepsilon,$$

*where the random variables* $X_t^{a,\psi}$ *are defined by* (2) *with* $a = a(\delta, K) = \delta K / \sqrt{\log(1 + N_\Pi(\delta/2))}$.

(ii) *For given DFP* $\Pi$ *of* $\Psi$ *and any constants* $K, L > 0$,

$$E^* \sup_{t \in [0,\tau]} \sup_{\psi, \phi \in \Psi} |X_t^\psi - X_t^\phi| 1_{\{\|\mathcal{W}\|_{\Pi,\tau} \leq K, |\overline{W}|^2 * \nu_\tau \leq L\}}$$

$$\lesssim K \int_0^{\Delta_\Pi} \sqrt{\log(1 + N_\Pi(\varepsilon))} \, d\varepsilon + \frac{L}{\Delta_\Pi K},$$

*where the random variables* $X_t^\psi$ *are defined by* (1).

Case A was already proved by [3] (Theorem 2.5), and the proof of Case B will be given in the Appendix.

By using the above result, we present a sufficient condition for the processes $\psi \rightsquigarrow X_\tau^\psi$ and $(t, \psi) \rightsquigarrow X_t^\psi$ to have paths which are bounded, almost surely. Actually, this property is trivial in Case A. On the other hand, in Case B, the result below gives the starting point of weak convergence theory in $\ell^\infty$-spaces.

THEOREM 2.6. *Consider Case* B. *For a given DFP $\Pi$ of $\Psi$, suppose that $\|\mathcal{W}\|_{\Pi,\tau} < \infty$ almost surely, and that $\int_0^{\Delta_\Pi} \sqrt{\log N_\Pi(\varepsilon)}\,d\varepsilon < \infty$. Then, the process $\psi \rightsquigarrow X_\tau^\psi$ takes values in $\ell^\infty(\Psi)$, almost surely. Furthermore, if the stopping time is a fixed time $\tau = t_0$, then the process $(t, \psi) \rightsquigarrow X_t^\psi$ takes values in $\ell^\infty([0, t_0] \times \Psi)$, almost surely.*

The proof will be given in the Appendix.

Now let us address to the context of weak convergence. Let $(E, \mathcal{E})$ be a Blackwell space, $\Psi$ an arbitrary set and $\Pi$ be a DFP of $\Psi$. We consider a sequence of the objects appearing above. That is, for every $n \in \mathbb{N}$, let $\mu^n$, $\nu^n$, $\mathcal{W}^n = \{W^{n,\psi} : \psi \in \Psi\}$, $\overline{W^n}$, $\|\mathcal{W}^n\|_\Pi$ and $\tau^n$ be the same objects as above with the new suffix $n$ for the sequence. Recall the definitions of Case $A^n$ and Case $B^n$ given in Section 1. Notice that Case $A^n$ is the same as $(3.1) + (3.2)$ of [3]. As a consequence of Theorems 2.5 and 2.6, we have the following claim.

COROLLARY 2.7. *The same assertions as Theorems* 3.2 *and* 3.4 *of* [3] *hold not only in Case* $A^n$ *but also in Case* $B^n$.

REMARK. In Theorem 2.5, no metric is equipped for $\Psi$. However, as in [3], we can use the tightness criterion based on partitioning given by van der Vaart and Wellner [6] (Theorem 1.5.6) rather than the well-known stochastic equicontinuity criterion.

## APPENDIX: PROOFS

PROOF OF THEOREM 2.5(i) IN CASE B. Fix any $\delta, K > 0$; we may assume $\delta \in \mathbb{Q}$ without loss of generality. For every integer $p \geq 0$, we define $a_p$, $\pi_p\psi$ and $\Pi_p\psi$ as in [3]. Fix any integer $q \geq 1$. For every $p = 0, 1, \ldots, q$, define $A_p(\psi)$ and $B_p(\psi)$ as in [3]. We do *not* introduce the stopping time $\tau_q$, and (2.4) in [3] should be read replacing $\tau_q$ by $\tau$.

Here, consider the identity

$$W^\psi - W^{\pi_0\psi} = (W^\psi - W^{\pi_0\psi})B_0(\psi) + \cdots$$

given by the lines 25–28 on page 691 of [3]. We have

$$\sup_{t\in[0,\tau]}\sup_{\psi\in\Psi}|X_t^{a(\delta,K),\psi}-X_t^{a(\delta,K),\pi_0\psi}|\le(I_1)+(I_2)+(II)+(III),$$

where $(I_1)$, $(I_2)$ and $(III)$ are from [3], and where

$$(II)=\sup_{t\in[0,\tau]}\sup_{\psi\in\Psi}|(W^\psi-W^{\pi_q\psi})A_q(\psi)*(\mu-\nu)_t|.$$

Notice that, in Case B, we do *not* have the inequality "$(II)\le(II_1)+(II_2)$", in which $(II_1)$ and $(II_2)$ are from [3]. This point is the difference between [3] and the present work.

The bounds for the terms $(I_1),(I_2)$ and $(III)$ are obtained by exactly the same way as [3]. On the other hand, if $\Psi$ is finite, then the term $(II)$ disappears as $q\to\infty$, because each $\Psi(2^{-q}\delta;k)$ contains only one point for sufficiently large $q$. If $\Psi$ is countable, that is, not finite, take a sequence $\{\Psi^m\}$ of finite subsets of $\Psi$ such that $\Psi^m\uparrow\Psi$. The proof is complete. $\square$

PROOF OF THEOREM 2.5(ii) IN CASE B.   The same as Theorem 2.5(ii) of [3]. $\square$

PROOF OF THEOREM 2.6.   We may assume that $\Pi$ is a NFP without loss of generality (see Lemma 2.4 of [3]), and we will use (i) of the above theorem for $\delta=\Delta_\Pi$ [note $N_\Pi(\Delta_\Pi)=1$]. Notice that $X^\psi=X^{a,\psi}+\check{X}^{a,\psi}$ almost surely, for any $a>0$. For any $K>0$, set $a=a(K)=a(\Delta_\Pi,K)=\Delta_\Pi K/\sqrt{\log(1+N_\Pi(\Delta_\Pi/2))}$.

First, we have

$$\sup_{t\in[0,\tau]}\sup_{\psi\in\Psi}|\check{X}_t^{a(K),\psi}|\le\overline{W}1_{\{\overline{W}>a(K)\}}*\mu_\tau+\overline{W}1_{\{\overline{W}>a(K)\}}*\nu_\tau.$$

Since $t\rightsquigarrow(\overline{W}^2\wedge\overline{W})*\nu_t$ is locally integrable, there exists an increasing sequence $\{T_m\}$ of stopping times such that, $E(\overline{W}1_{\{\overline{W}>a(K)\}}*\nu_{T_m\wedge\tau})<\infty$, thus, $\overline{W}1_{\{\overline{W}>a(K)\}}*\nu_{T_m\wedge\tau}<\infty$ almost surely. By letting $m\to\infty$, we have $\overline{W}1_{\{\overline{W}>a(K)\}}*\nu_\tau<\infty$ almost surely. Since the residual $\overline{W}1_{\{\overline{W}>a(K)\}}*(\mu-\nu)$ is a local martingale, we also have $\overline{W}1_{\{\overline{W}>a(K)\}}*\mu_\tau<\infty$ almost surely. Hence, it holds that $\sup_{t\in[0,\tau]}\sup_{\psi\in\Psi}|\check{X}_t^{a(K),\psi}|<\infty$ almost surely.

Next, by (i) of Theorem 2.5, we have

$$E\sup_{t\in[0,\tau]}\sup_{\psi,\phi\in\Psi}|X_t^{a(K),\psi}-X_t^{a(K),\phi}|1_{\{\|\mathcal{W}\|_{\Pi,\tau}\le K\}}$$
$$\lesssim K\int_0^{\Delta_\Pi}\sqrt{\log(1+N_\Pi(\varepsilon))}\,d\varepsilon<\infty.$$

Hence, we have

$$\sup_{t \in [0,\tau]} \sup_{\psi \in \Psi} |X_t^{a(K),\psi}| 1_{\{\|\mathcal{W}\|_{\Pi,\tau} \leq K\}} < \infty \qquad \text{almost surely.}$$

So it holds that

$$P\left(\sup_{t \in [0,\tau]} \sup_{\psi \in \Psi} |X_t^\psi| = \infty, \|\mathcal{W}\|_{\Pi,\tau} \leq K\right) = 0 \qquad \forall K > 0.$$

We therefore obtain

$$P\left(\sup_{t \in [0,\tau]} \sup_{\psi \in \Psi} |X_t^\psi| = \infty\right) \leq \sum_{K=1}^\infty P\left(\sup_{t \in [0,\tau]} \sup_{\psi \in \Psi} |X_t^\psi| = \infty, \|\mathcal{W}\|_{\Pi,\tau} \leq K\right) = 0.$$

This proves the assertions of the theorem. $\square$

**Acknowledgments.** I thank the Editor and the referees for their comments and suggestions which improved the presentation of the paper.

Institute of Statistical Mathematics
4-6-7 Minami-Azabu, Minato-ku
Tokyo 106-8569
Japan
E-mail: nisiyama@ism.ac.jp